\documentclass{article}

\newcommand{\be}{\begin{equation}}
\newcommand{\ee}{\end{equation}}
\newcommand{\bea}{\begin{eqnarray}}
\newcommand{\eea}{\end{eqnarray}}
\newcommand{\barray}{\begin{array}}
\newcommand{\earray}{\end{array}}
\newcommand{\pa}{\partial}
\newcommand{\nn}{\nonumber}
\newcommand{\bitem}{\begin{itemize}}
\newcommand{\eitem}{\end{itemize}}
\newtheorem{teo}{Theorem}[section]
\newcommand{\bt}{\begin{teo}}
\newcommand{\et}{\end{teo}}
\newtheorem{Def}{Definition}[section]
\newcommand{\bd}{\begin{Def}}
\newcommand{\ed}{\end{Def}}
\newtheorem{lem}{Lemma}[section]
\newcommand{\bl}{\begin{lem}}
\newcommand{\el}{\end{lem}}
\newtheorem{prop}{Proposition}[section]
\newcommand{\bp}{\begin{prop}}
\newcommand{\ep}{\end{prop}}
\newtheorem{cor}{Corollary}[section]
\newcommand{\bc}{\begin{cor}}
\newcommand{\ec}{\end{cor}}
\newtheorem{ex}{Example}[section]
\newcommand{\bex}{\begin{ex}}
\newcommand{\eex}{\end{ex}}
\newtheorem{rem}{Remark}[section]
\newcommand{\br}{\begin{rem}}
\newcommand{\er}{\end{rem}}

\catcode `\@=11
\@addtoreset{equation}{section}

\begin{document}

\begin{center}
{\Large \textbf{Compatible nonlocal Poisson brackets \\
of hydrodynamic
type, and integrable \\ hierarchies
related to them\footnote{This work was supported by
the Alexander von Humboldt Foundation (Germany),
the Russian Foundation for Basic Research
(grant No. 99--01--00010) and the INTAS
(grant No. 99--1782).}}}
\end{center}

\bigskip
\bigskip

\centerline{\large {O. I. Mokhov}}
\bigskip
\medskip

\section{Introduction. Basic definitions} \label{vved}

In the present work, the integrable bi-Hamiltonian
hierarchies related to compatible nonlocal Poisson
brackets of hydrodynamic type are effectively
constructed.
For achieving this aim, first of all, the problem on
the canonical form of a special type for
compatible nonlocal Poisson brackets of hydrodynamic type
is solved. The compatible pairs of nonlocal
Poisson brackets of hydrodynamic type
have a more simple description in special coordinates
in which the metrics corresponding to these brackets
are diagonal (see \cite{1}, \cite{2}), but
for an effective construction of the hierarchies
we need a different approach developed in this paper.
For compatible local Poisson brackets of
hydrodynamic type (the Dubrovin--Novikov brackets \cite{3}),
the corresponding integrable bi-Hamiltonian
hierarchies were constructed by the present author in
the papers \cite{4}, \cite{5}, and for compatible
nonlocal Mokhov--Ferapontov brackets \cite{6}
generated by metrics of constant Riemannian curvature,
the bi-Hamiltonian hierarchies were constructed in \cite{7}.

\subsection{Local Poisson brackets of hydrodynamic type}

An arbitrary local homogeneous
first-order Poisson bracket, that is,
a Poisson bracket of the form
\be
\{ u^i (x), u^j (y) \} =
g^{ij} (u(x))\, \delta_x (x-y) + b^{ij}_k (u(x)) \, u^k_x \,
\delta (x-y), \label{(1.1)}
\ee
where $u^1,...,u^N$ are local coordinates
on a certain given smooth $N$-dimensional
manifold $M$, is called {\it a local Poisson bracket of
hydrodynamic type} or {Dubrovin--Novikov bracket} \cite{3}.
Here $u^i(x),\ 1 \leq i \leq N,$ are functions (fields)
of single independent variable $x$, the coefficients
$g^{ij}(u)$ and $b^{ij}_k (u)$ of bracket (\ref{(1.1)})
are smooth functions of local coordinates.

In other words,
for arbitrary functionals $I[u]$ and $J[u]$
on the space of fields $u^i(x), \ 1 \leq i \leq N,$ a bracket of the form
\be
\{ I,J \} = \int {\delta I \over \delta u^i(x)}
\biggl ( g^{ij}(u(x)) {d \over dx} + b^{ij}_k (u(x))\, u^k_x \biggr )
{\delta J \over \delta u^j(x)} dx
\label{lok}
\ee
is defined and it is required that this bracket is a Poisson bracket,
that is, it is skew-symmetric:
\be
\{ I, J \} = - \{ J, I \}, \label{skew}
\ee
and satisfies the Jacobi identity
\be
\{ \{ I, J \}, K \} + \{ \{ J, K \}, I \} + \{ \{ K, I \}, J \} =0
\label{jacobi}
\ee
for arbitrary functionals $I[u]$, $J[u]$, and $K[u]$.

A local bracket (\ref{lok}) is called
{\it nondegenerate} if
$\det (g^{ij} (u)) \not\equiv 0$.
For the general nondegenerate brackets (\ref{lok}),
Dubrovin and Novikov proved the following important theorem.

\bt [Dubrovin, Novikov \cite{3}] \label{dn}
If $\det (g^{ij} (u)) \not\equiv 0$, then bracket (\ref{lok})
is a Poisson bracket, that is, it is skew-symmetric and satisfies
the Jacobi identity, if and only if
\bitem
\item [(1)] $g^{ij} (u)$ is an arbitrary flat
pseudo-Riemannian contravariant metric (a metric of
zero Riemannian curvature),

\item [(2)] $b^{ij}_k (u) = - g^{is} (u) \Gamma ^j_{sk} (u),$ where
$\Gamma^j_{sk} (u)$ is the Riemannian connection generated by
the contravariant metric $g^{ij} (u)$
(the Levi--Civita connection).
\eitem
\et

Consequently, for any local nondegenerate Poisson bracket of
hydrodynamic type, there always exist local coordinates
$v^1,...,v^N$ (flat coordinates of the metric $g^{ij}(u)$) in which
all the coefficients of the bracket are constant:
$$
\widetilde g^{ij} (v)
= \eta^{ij} = {\rm \ const}, \ \
\widetilde \Gamma^i_{jk} (v) = 0, \ \
\widetilde b^{ij}_k (v) =0,
$$
that is, the bracket has the form
\be
\{ I,J \} = \int {\delta I \over \delta v^i(x)}
 \eta^{ij} {d \over dx}
{\delta J \over \delta v^j(x)} dx,
\label{(1.6)}
\ee
where $(\eta^{ij})$ is a nondegenerate symmetric constant matrix:
$$
\eta^{ij} = \eta^{ji}, \ \ \eta^{ij} = {\rm const},
\  \ \det \, (\eta^{ij}) \neq 0.
$$

The local Poisson brackets of hydrodynamic type (\ref{(1.1)})
were introduced and studied by Dubrovin and Novikov
in \cite{3}. In this paper, they proposed
a general local Hamiltonian approach (this approach
corresponds to the local brackets of form (\ref{(1.1)}))
to the so-called
{\it homogeneous systems of hydrodynamic type}, that is, to
evolutionary quasilinear systems of first-order partial differential
equations
\be
u^i_t = V^i_j (u)\, u^j_x.
\label{(1.5)}
\ee

This Hamiltonian approach was motivated by the study
of the equations of Euler hydrodynamics and the
Whitham averaging equations describing the evolution of
slowly modulated multiphase solutions of partial
differential equations (see \cite{8}).
In \cite{9}, \cite{10} Tsarev constructed the
theory of integrating the class of diagonalizable
Hamiltonian (and also semi-Hamiltonian)
homogeneous systems of hydrodynamic type.

\subsection{Nonlocal Poisson brackets of
hydrodynamic type}

Nonlocal Poisson brackets of hydrodynamic type
(the Mokhov--Ferapontov brackets) were
introduced and studied in the work of the present author
and Ferapontov \cite{6}
(see also \cite{11}--\cite{13}).
They have the following form:

\be
\{ I,J \} = \int {\delta I \over \delta u^i(x)}
\left ( g^{ij}(u(x)) {d \over dx} + b^{ij}_k (u(x))\, u^k_x
+ K u^i_x \left ( {d \over dx} \right )^{-1} u^j_x \right )
{\delta J \over \delta u^j(x)} dx,
\label{nonl}
\ee
where $K$ is an arbitrary constant.

A bracket of form (\ref{nonl}) is called {\it nondegenerate}
if
$\det (g^{ij} (u)) \not\equiv 0$.

\bt [\cite{6}] \label{mofer}
If $\det (g^{ij} (u)) \not\equiv 0$, then bracket (\ref{nonl})
is a Poisson bracket, that is, it is skew-symmetric and
satisfies the Jacobi identity, if and only if
\bitem
\item [(1)] $g^{ij} (u)$ is an arbitrary pseudo-Riemannian
contravariant metric of
constant Riemannian curvature $K$,

\item [(2)] $b^{ij}_k (u) = - g^{is} (u) \Gamma ^j_{sk} (u),$ where
$\Gamma^j_{sk} (u)$ is the Riemannian connection
generated by the contravariant metric $g^{ij} (u)$
(the Levi--Civita connection).
\eitem
\et

In \cite{6}
Ferapontov introduced and studied more general
nonlocal Poisson brackets of hydrodynamic type
(the Ferapontov brackets), namely,
the Poisson brackets of the form
\bea
&& \
\{ I,J \} = \int   {\delta I \over \delta u^i(x) }
\left ( g^{ij}(u(x)) {d \over dx} +
b^{ij}_k (u(x))\, u^k_x + \right.
\nn\\
&& \
\left. \sum_{\alpha =1}^L
\varepsilon_{\alpha} (w^{\alpha})^i_k (u (x))
u^k_x \left ( {d \over dx} \right )^{-1}
(w^{\alpha})^j_s (u (x)) u^s_x \right )
{\delta J \over \delta u^j(x)} dx,
\ \ \det (g^{ij} (u)) \not\equiv 0.
\label{nonl2}
\eea

\bt [\cite{11}] \label{ferap}
Bracket (\ref{nonl2}) is a Poisson bracket,
that is, it is skew-symmetric and satisfies the
Jacobi identity, if and only if
\bitem
\item [(1)] $b^{ij}_k (u) = - g^{is} (u) \Gamma ^j_{sk} (u),$
where
$\Gamma^j_{sk} (u)$ is the Riemannian connection
generated by the contravariant metric $g^{ij} (u)$
(the Levi--Civita connection),

\item [(2)] the metric $g^{ij} (u)$ and the set of the affinors
$(w^{\alpha})^i_j (u)$ satisfies relations

\be
g_{ik} (u) (w^{\alpha})^k_j (u) =
g_{jk} (u) (w^{\alpha})^k_i (u),
\ \ \ \alpha = 1,...,L,  \label{peter1}
\ee
\be
\nabla_k (w^{\alpha})^i_j (u) =
\nabla_j (w^{\alpha})^i_k (u),
\ \ \ \alpha = 1,...,L,   \label{peter2}
\ee
\be
R^{ij}_{kl} (u) =  \sum_{\alpha =1}^L
\varepsilon_{\alpha}
\left ( (w^{\alpha})^i_l (u) (w^{\alpha})^j_k (u)
- (w^{\alpha})^j_l (u) (w^{\alpha})^i_k (u) \right ). \label{gauss}
\ee
In addition, the family of the affinors $w^{\alpha} (u)$
is commutative: $[w^{\alpha}, w^{\beta}] =0.$
\eitem
\et

Let us write out all the relations on the coefficients
of the nonlocal Poisson bracket
(\ref{nonl2}) in a convenient form for further
repeated use.

\bl
Bracket (\ref{nonl2}) is a Poisson bracket
if and only if its coefficients satisfy the
relations
\be
g^{ij} = g^{ji}, \label{01}
\ee
\be
{\pa g^{ij} \over \pa u^k} = b^{ij}_k + b^{ji}_k, \label{02}
\ee
\be
g^{is} b^{jk}_s = g^{js} b^{ik}_s, \label{03}
\ee
\be
g^{is} (w^{\alpha})^j_s = g^{js} (w^{\alpha})^i_s, \label{04}
\ee
\be
(w^{\alpha})^i_s (w^{\beta})^s_j =
(w^{\beta})^i_s (w^{\alpha})^s_j, \label{05}
\ee
\be
g^{is} g^{jr} {\pa (w^{\alpha})^k_r \over \pa u^s}
- g^{jr} b^{ik}_s (w^{\alpha})^s_r  =
g^{js} g^{ir} {\pa (w^{\alpha})^k_r \over \pa u^s}
- g^{ir} b^{jk}_s (w^{\alpha})^s_r, \label{06}
\ee
\be
g^{is} \left ( {\pa b^{jk}_r \over \pa u^s}
- {\pa b^{jk}_s \over \pa u^r} \right )
+ b^{ik}_s b^{sj}_r - b^{ij}_s b^{sk}_r =
\sum_{\alpha = 1}^L \varepsilon_{\alpha} g^{is}
\left ( (w^{\alpha})^j_s (w^{\alpha})^k_r -
 (w^{\alpha})^j_r (w^{\alpha})^k_s \right ). \label{07}
\ee
\el

\subsection{Compatible Poisson brackets}

In \cite{14} Magri proposed
a bi-Hamiltonian approach to the integration of nonlinear
systems. This approach demonstrated
that integrability is closely related to the
bi-Hamiltonian property, that is, the property
of a system to have two compatible Hamiltonian
representations.

\bd [Magri \cite{14}] \label{dm}
{\rm Two Poisson brackets $\{ \, \cdot \, , \, \cdot \, \}_1$
and $\{ \, \cdot \, , \, \cdot \, \}_2$ are called {\it compatible} if
an arbitrary linear combination of these Poisson brackets
\be
\{ \, \cdot \, , \, \cdot \, \} =
\lambda_1 \, \{ \, \cdot \, , \, \cdot \, \}_1 +
\lambda_2 \, \{ \, \cdot \, , \, \cdot \, \}_2, \label{magri}
\ee
where $\lambda_1$ and $\lambda_2$ are arbitrary constants,
is also a Poisson bracket.
In this case, we shall also say that the brackets
$\{ \, \cdot \, , \, \cdot \, \}_1$ and
$\{ \, \cdot \, , \, \cdot \, \}_2$
form {\it a pencil of Poisson brackets}.}
\ed

 As was shown by Magri in \cite{14},
compatible Poisson brackets generate integrable hierarchies of
systems of differential equations.
In particular, for a system, the bi-Hamiltonian property
generates recurrent relations for the conservation laws of
this system.

Here the integrable hierarchies related to compatible
nonlocal Poisson brackets of hydrodynamic type are constructed.

\section{Pencil of nonlocal Poisson brackets\\
of hydrodynamic type}

Let us describe all nonlocal Poisson brackets (\ref{nonl2})
compatible with the constant nondegenerate Poisson bracket
of hydrodynamic type
\be
\ \{ I, J \}_2 = \int {\delta I \over \delta u^i(x)}
\eta^{ij} {d \over dx}
{\delta J \over \delta u^j(x)} dx, \label{co}
\ee
where $(\eta^{ij})$ is an arbitrary nondegenerate
symmetric constant matrix: $\det (\eta^{ij}) \neq 0,$
$\eta^{ij} = \eta^{ji},$ $\eta^{ij} = {\rm const},$
that is,
let us classify all the following pencils
of nonlocal Poisson brackets:
\be
\ \{ I, J \}_{\lambda} = \{ I, J \}_1 + \lambda \{ I, J \}_2,
\ee
where $\{ I, J \}_1$ is a Poisson bracket of form (\ref{nonl2}).

\bl
The Poisson brackets (\ref{co}) and (\ref{nonl2})
are compatible if and only if the following relations are
satisfied:
\be
\eta^{is} b^{jk}_s = \eta^{js} b^{ik}_s,  \label{1}
\ee
\be
\eta^{is} (w^{\alpha})^j_s = \eta^{js} (w^{\alpha})^i_s,
\label{2}
\ee
\be
{\pa (w^{\alpha})^i_j \over \pa u^k}
= {\pa (w^{\alpha})^i_k \over \pa u^j}, \label{3}
\ee
\be
{\pa b^{jk}_r \over \pa u^s}
- {\pa b^{jk}_s \over \pa u^r}
=
\sum_{\alpha = 1}^L \varepsilon_{\alpha}
\left ( (w^{\alpha})^j_s (w^{\alpha})^k_r -
 (w^{\alpha})^j_r (w^{\alpha})^k_s \right ). \label{5}
\ee
\el

It is important to note that the relation
\be
\eta^{jr} b^{ik}_s (w^{\alpha})^s_r  =
\eta^{ir} b^{jk}_s (w^{\alpha})^s_r,  \label{4}
\ee
derived from (\ref{06}) as one of the compatibility
conditions for Poisson brackets (\ref{co}) and (\ref{nonl2})
follows from relations (\ref{1})--(\ref{3})
for every Poisson bracket (\ref{nonl2}).
Actually, from (\ref{2}) we get
\be
\eta^{jr} b^{ik}_s (w^{\alpha})^s_r  =
\eta^{sr} b^{ik}_s (w^{\alpha})^j_r,
\ee
and from relation (\ref{1}) we also have
\be
\eta^{sr} b^{ik}_s (w^{\alpha})^j_r=
\eta^{si} b^{rk}_s (w^{\alpha})^j_r.
\ee
Consequently relation (\ref{4}) is reduced to the
relation
\be
\eta^{si} b^{rk}_s (w^{\alpha})^j_r =
\eta^{ir} b^{jk}_s (w^{\alpha})^s_r,
\ee
that is,
\be
b^{rk}_s (w^{\alpha})^j_r =
b^{jk}_r (w^{\alpha})^r_s. \label{bw}
\ee
Let us prove that relation (\ref{bw})
is satisfied for any Poisson bracket (\ref{nonl2})
for which relation (\ref{3}) is valid.
In fact, in this case relation (\ref{06})
for the Poisson bracket (\ref{nonl2}) takes the form
\be
g^{jr} b^{ik}_s (w^{\alpha})^s_r  =
g^{ir} b^{jk}_s (w^{\alpha})^s_r. \label{gbw}
\ee
From relation (\ref{04}) we get
\be
g^{jr} b^{ik}_s (w^{\alpha})^s_r  =
g^{sr} b^{ik}_s (w^{\alpha})^j_r,
\ee
and from relation (\ref{03}) we have
\be
g^{sr} b^{ik}_s (w^{\alpha})^j_r =
g^{si} b^{rk}_s (w^{\alpha})^j_r,
\ee
that is, relation (\ref{gbw}) is reduced to
the relation
\be
g^{si} b^{rk}_s (w^{\alpha})^j_r =
g^{ir} b^{jk}_s (w^{\alpha})^s_r,
\ee
which is equivalent to relation (\ref{bw}).

\section{Canonical form of compatible
pairs of brackets}

\bt    \label{t1}
An arbitrary nonlocal Poisson bracket $\{ I, J \}_1$ of form
(\ref{nonl2}) is compatible with the constant Poisson bracket
(\ref{co}) if and only if it has the form
\bea
&& \
\{ I, J \}_1 = \int   {\delta I \over \delta u^i(x) }
\left (
\left [ \eta^{is}
{\pa F^j \over \pa u^s} +
\eta^{js}
{\pa F^i \over \pa u^s} -
\eta^{is} \eta^{jk}
\sum_{\alpha = 1}^L
\varepsilon_{\alpha} {\pa \psi^{\alpha} \over \pa u^s}
{\pa \psi^{\alpha} \over \pa u^k} \right ]
 {d \over dx} + \right. \nn\\
&& \
\left [ \eta^{is} {\pa^2 F^j \over \pa u^s \pa u^k}
- \eta^{is} \eta^{jp} \sum_{\alpha = 1}^L
\varepsilon_{\alpha} {\pa^2 \psi^{\alpha} \over
\pa u^s \pa u^k} {\pa \psi^{\alpha} \over \pa u^p} \right  ]
\, u^k_x + \nn\\
&& \
\left. \eta^{ip} \eta^{jr}
\sum_{\alpha =1}^L
\varepsilon_{\alpha} {\pa^2 \psi^{\alpha} \over \pa u^p \pa u^k}
u^k_x \left ( {d \over dx} \right )^{-1}
{\pa^2 \psi^{\alpha} \over \pa u^r \pa u^s} u^s_x \right )
{\delta J \over \delta u^j(x)} dx, \label{nonloc3}
\eea
where $F^i (u),$ $1 \leq i \leq N,$ and $\psi^{\alpha} (u),$
$1 \leq \alpha \leq L,$ are smooth functions
defined in a certain domain of local coordinates.
\et

It follows immediately from relation (\ref{3})
that there locally exist functions $(\varphi^{\alpha})^i (u),$
$1 \leq i \leq N,$ $1 \leq \alpha \leq L,$ such that
\be
(w^{\alpha})^i_j =
{\pa (\varphi^{\alpha})^i \over \pa u^j}.
\ee

Then relation (\ref{5}) takes the form
\be
{\pa b^{jk}_r \over \pa u^s}
- {\pa b^{jk}_s \over \pa u^r}
=
\sum_{\alpha = 1}^L \varepsilon_{\alpha}
\left ( {\pa (\varphi^{\alpha})^j \over \pa u^s}
{(\pa \varphi^{\alpha})^k \over \pa u^r } -
 {\pa (\varphi^{\alpha})^j \over \pa u^r}
{\pa (\varphi^{\alpha})^k \over \pa u^s} \right ). \label{5a}
\ee

Let us introduce the function
\be
A^{ij}_k (u) =
b^{ij}_k (u) - \sum_{\alpha = 1}^L \varepsilon_{\alpha}
(\varphi^{\alpha})^i {\pa (\varphi^{\alpha})^j \over \pa u^k}.
\ee
Then using (\ref{5a}) we get
\bea
&&
{\pa A^{ij}_k \over \pa u^l} =
{\pa b^{ij}_k \over \pa u^l} - \sum_{\alpha = 1}^L
\varepsilon_{\alpha} {\pa (\varphi^{\alpha})^i \over \pa u^l}
{\pa (\varphi^{\alpha})^j \over \pa u^k} -
\sum_{\alpha = 1}^L
\varepsilon_{\alpha} (\varphi^{\alpha})^i
{\pa^2 (\varphi^{\alpha})^j \over \pa u^k \pa u^l}=\nn\\
&&
{\pa b^{ij}_l \over \pa u^k} - \sum_{\alpha = 1}^L
\varepsilon_{\alpha} {\pa (\varphi^{\alpha})^i \over \pa u^k}
{\pa (\varphi^{\alpha})^j \over \pa u^l} -
\sum_{\alpha = 1}^L
\varepsilon_{\alpha} (\varphi^{\alpha})^i
{\pa^2 (\varphi^{\alpha})^j \over \pa u^k \pa u^l} =
{\pa A^{ij}_l \over \pa u^k}.
\eea
Consequently there exist functions $P^{ij} (u),$
$1 \leq i, j \leq N,$ such that
\be
A^{ij}_k (u) = {\pa P^{ij} \over \pa u^k}.
\ee
Thus,
\be
b^{ij}_k (u) = {\pa P^{ij} \over \pa u^k} +
\sum_{\alpha = 1}^L \varepsilon_{\alpha}
(\varphi^{\alpha})^i {\pa (\varphi^{\alpha})^j \over \pa u^k}.
\ee
From relation (\ref{02}) for the Poisson bracket
$\{ I, J \}_1$ we get
\be
{\pa g^{ij} \over \pa u^k} = b^{ij}_k + b^{ji}_k =
{\pa P^{ij} \over \pa u^k} +
{\pa P^{ji} \over \pa u^k} +
\sum_{\alpha = 1}^L \varepsilon_{\alpha}
(\varphi^{\alpha})^i {\pa (\varphi^{\alpha})^j \over \pa u^k}
+ \sum_{\alpha = 1}^L \varepsilon_{\alpha}
(\varphi^{\alpha})^j {\pa (\varphi^{\alpha})^i \over \pa u^k},
\ee
that is, using (\ref{01}) we have
\be
g^{ij} =
P^{ij} + P^{ji} +
\sum_{\alpha = 1}^L \varepsilon_{\alpha}
(\varphi^{\alpha})^i (\varphi^{\alpha})^j + c^{ij} +c^{ji},
\ee
where $(c^{ij})$ is an arbitrary constant matrix, $c^{ij} =
{\rm const}.$
Defining the function $R^{ij} (u)$ by the formula
\be
R^{ij} = P^{ij} + c^{ij},
\ee
we get the proof of the Liouville property for
the Poisson bracket $\{ I, J \}_1$ in the
considered local coordinates:
\be
g^{ij} =
R^{ij} + R^{ji} +
\sum_{\alpha = 1}^L \varepsilon_{\alpha}
(\varphi^{\alpha})^i (\varphi^{\alpha})^j,
\ee
\be
b^{ij}_k = {\pa R^{ij} \over \pa u^k} +
\sum_{\alpha = 1}^L \varepsilon_{\alpha}
(\varphi^{\alpha})^i {\pa (\varphi^{\alpha})^j \over \pa u^k}
\ee
(in the next section, see about the Liouville property
more in detail).

It follows from relation (\ref{2}) that
\be
\eta^{is} {\pa ( \varphi^{\alpha})^j \over \pa u^s} =
\eta^{js} {\pa ( \varphi^{\alpha})^i \over \pa u^s},
\ee
that is,
\be
{\pa \left ( \eta_{rj} ( \varphi^{\alpha})^j
\right ) \over \pa u^p} =
{\pa \left ( \eta_{pj} ( \varphi^{\alpha})^j
\right ) \over \pa u^r}.
\ee
Thus there exist functions $\psi^{\alpha} (u),$
$1 \leq \alpha \leq L,$ such that
\be
\eta_{rs} (\varphi^{\alpha})^s =
{\pa \psi^{\alpha} \over \pa u^r}
\ee
or
\be
 (\varphi^{\alpha})^s =
\eta^{sr} {\pa \psi^{\alpha} \over \pa u^r}.
\ee

It follows from relation (\ref{1}) that
\be
\eta^{is} \left (
{\pa R^{jk} \over \pa u^s} + \eta^{jr} \eta^{kp}
\sum_{\alpha = 1}^L
\varepsilon_{\alpha} {\pa \psi^{\alpha} \over \pa u^r}
{\pa^2 \psi^{\alpha} \over \pa u^p \pa u^s} \right ) =
\eta^{js} \left (
{\pa R^{ik} \over \pa u^s} + \eta^{ir} \eta^{kp}
\sum_{\alpha = 1}^L
\varepsilon_{\alpha} {\pa \psi^{\alpha} \over \pa u^r}
{\pa^2 \psi^{\alpha} \over \pa u^p \pa u^s} \right )
\ee
or
\be
{\pa ( \eta_{lj} R^{jk}) \over \pa u^s} + \eta^{kp}
\sum_{\alpha = 1}^L
\varepsilon_{\alpha} {\pa \psi^{\alpha} \over \pa u^l}
{\pa^2 \psi^{\alpha} \over \pa u^p \pa u^s}  =
{\pa ( \eta_{sj} R^{jk}) \over \pa u^l} + \eta^{kp}
\sum_{\alpha = 1}^L
\varepsilon_{\alpha} {\pa \psi^{\alpha} \over \pa u^s}
{\pa^2 \psi^{\alpha} \over \pa u^p \pa u^l}.
\ee
Thus,
\be
{\pa \left ( \eta_{lj} R^{jk} + \eta^{kp}
\sum_{\alpha = 1}^L
\varepsilon_{\alpha} {\pa \psi^{\alpha} \over \pa u^l}
{\pa \psi^{\alpha} \over \pa u^p} \right )
\over \pa u^s} =
{\pa \left ( \eta_{sj} R^{jk} + \eta^{kp}
\sum_{\alpha = 1}^L
\varepsilon_{\alpha} {\pa \psi^{\alpha} \over \pa u^s}
{\pa \psi^{\alpha} \over \pa u^p } \right )
\over \pa u^l}.
\ee
Consequently there exist functions $F^j (u)$
such that
\be
\eta_{lj} R^{jk} + \eta^{kp}
\sum_{\alpha = 1}^L
\varepsilon_{\alpha} {\pa \psi^{\alpha} \over \pa u^l}
{\pa \psi^{\alpha} \over \pa u^p} = {\pa F^k \over \pa u^l},
\ee
that is, we get
\be
R^{ij} = \eta^{is}
{\pa F^j \over \pa u^s} -
\eta^{is} \eta^{jk}
\sum_{\alpha = 1}^L
\varepsilon_{\alpha} {\pa \psi^{\alpha} \over \pa u^s}
{\pa \psi^{\alpha} \over \pa u^k},
\ee
\be
g^{ij} =
\eta^{is}
{\pa F^j \over \pa u^s} +
\eta^{js}
{\pa F^i \over \pa u^s} -
\eta^{is} \eta^{jk}
\sum_{\alpha = 1}^L
\varepsilon_{\alpha} {\pa \psi^{\alpha} \over \pa u^s}
{\pa \psi^{\alpha} \over \pa u^k},
\ee
\be
b^{ij}_k = \eta^{is} {\pa^2 F^j \over \pa u^s \pa u^k}
- \eta^{is} \eta^{jp} \sum_{\alpha = 1}^L
\varepsilon_{\alpha} {\pa^2 \psi^{\alpha} \over
\pa u^s \pa u^k} {\pa \psi^{\alpha} \over \pa u^p}.
\ee

Here it is easy to check that for the Poisson bracket (\ref{nonloc3})
all the relations of compatibility (\ref{1})--(\ref{5})
are satisfied.

\section{Integrable equations for canonical \\
compatible pair of brackets}

\bt \label{liu}
Nonlocal bracket (\ref{nonloc3}) is a Poisson bracket
if and only if the following relations are satisfied:
\be
{\pa^2 Q_1 \over \pa u^i \pa u^s} \eta^{sp}
{\pa^2 Q_2 \over \pa u^p \pa u^j} =
{\pa^2 Q_2 \over \pa u^i \pa u^s} \eta^{sp}
{\pa^2 Q_1 \over \pa u^p \pa u^j},  \label{ass1}
\ee
\be
g^{is}_{{\rm  can}} \eta^{jr} {\pa^2 Q \over \pa u^r \pa u^s}
= g^{js}_{{\rm can}} \eta^{ir} {\pa^2 Q \over \pa u^r \pa u^s},
\label{ass2}
\ee
where $g^{ij}_{{\rm can}} (u)$ is the metric
of bracket (\ref{nonloc3}):
\be
g^{ij}_{{\rm can}} (u) =
\eta^{is}
{\pa F^j \over \pa u^s} +
\eta^{js}
{\pa F^i \over \pa u^s} -
\eta^{is} \eta^{jk}
\sum_{\alpha = 1}^L
\varepsilon_{\alpha} {\pa \psi^{\alpha} \over \pa u^s}
{\pa \psi^{\alpha} \over \pa u^k},
\ee
the functions $Q (u),$ $Q_1 (u),$ and $Q_2 (u)$ are
arbitrary from the functions
$F^i (u),$ $1 \leq i \leq N,$ and $\psi^{\alpha} (u),$
$1 \leq \alpha \leq L.$
\et

The nonlinear system (\ref{ass1}), (\ref{ass2}) is integrable
by the method of inverse scattering problem.
The procedure of the integration for the system (\ref{ass1}),
(\ref{ass2}) will be published in our next paper.

\section{Liouville and special Liouville coordinates}

Local coordinates $u = (u^1,..., u^N)$ are called
{\it Liouville} for an arbitrary Poisson bracket
$\{ I, J \}$ if the functions (the fields) $u^i (x)$
are densities of integrals in involution with
respect to this bracket, that is,
\be
\ \{ U^i, U^j \} = 0, \ \ \ 1 \leq i, j \leq N,
\ee
where $U^i = \int u^i (x) dx,$ $1 \leq i \leq N.$
In this case the Poisson bracket is also
called {\it Liouville} in these coordinates.
Liouville coordinates naturally arise and
play an essential role in the Dubrovin--Novikov
procedure of averaging of Hamiltonian equations \cite{3}.
Physical coordinates derived by averaging
of densities of participating in the Dubrovin--Novikov
procedure $N$
involutive local integrals of an initial
Hamiltonian system are always Liouville for
corresponding averaged bracket.
This property was a motivation for
the definition of Liouville coordinates
for local Poisson brackets of hydrodynamic type in \cite{3}.
For general nonlocal Poisson brackets of hydrodynamic type
(\ref{nonl2}),
Liouville coordinates were introduced in \cite{15}.

A nonlocal Poisson bracket of hydrodynamic type
(\ref{nonl2}) is Liouville in the local coordinates
$u = (u^1,...,u^N)$ if and only if
there exist functions
$(\varphi^{\alpha})^i (u),$ $1 \leq i \leq L$, and
a matrix function
$\Phi^{ij} (u)$ such that the bracket has the following form
(see \cite{15}, where namely this characteristic is taken
for the definition of the Liouville property of bracket
(\ref{nonl2})):
\bea
&& \
\{ I, J \}_1 = \int   {\delta I \over \delta u^i(x) }
\left (
\left [ \Phi^{ij} (u) + \Phi^{ji} (u) -
\sum_{\alpha = 1}^L
\varepsilon_{\alpha} (\varphi^{\alpha})^i
(\varphi^{\alpha})^j \right ]
 {d \over dx} + \right. \nn\\
&& \
\left [ {\pa \Phi^{ij} \over \pa u^k}
- \sum_{\alpha = 1}^L
\varepsilon_{\alpha} {\pa (\varphi^{\alpha})^i \over
\pa u^k} \varphi^{\alpha})^j \right  ]
\, u^k_x + \nn\\
&& \
\left.
\sum_{\alpha =1}^L
\varepsilon_{\alpha} {\pa (\varphi^{\alpha})^i \over \pa u^k}
u^k_x \left ( {d \over dx} \right )^{-1}
{\pa (\varphi^{\alpha})^j \over \pa u^s} u^s_x \right )
{\delta J \over \delta u^j(x)} dx. \label{liuv0}
\eea

From theorem \ref{t1}, it follows

\bt
Flat coordinates of an arbitrary nondegenerate
local Poisson bracket of hydrodynamic type $\{ I, J \}_2$
are always Liouville for any nonlocal Poisson bracket
$\{ I, J \}_1$ (\ref{nonl2}) compatible with $\{ I, J \}_2$.
Moreover, in addition the corresponding Liouville
function
$\Phi^{ij} (u)$ and the functions
$(\varphi^{\alpha})^i (u)$ always have the
special form
\be
\Phi^{ij} (u) = \eta^{is} {\pa F^j \over \pa u^s},\ \ \ \
(\varphi^{\alpha})^i (u) = \eta^{is} {\pa \psi^{\alpha} \over \pa u^s}.
\ee
\et

Local coordinates $u = (u^1,...,u^N)$ are called
{\it special Liouville coordinates} \cite{16},
\cite{17}
for an arbitrary Poisson bracket
$\{ I, J \}$ if there exists a nonzero
constant symmetric matrix $(\eta_{ij})$ such that
the functions (the fields) $u^i (x),$ $1 \leq i \leq N,$
and $\eta_{ij} u^i (x) u^j (x)$ are densities of
integrals in involution with respect to this bracket, that is,
\be
\ \{ U^i, U^j \} = 0, \ \ \ 1 \leq i, j \leq N + 1,
\ee
where $U^i = \int u^i (x) dx,$ $1 \leq i \leq N,$
$U^{N+1} = \int \eta_{ij} u^i (x) u^j (x) dx$.
In this case the Poisson bracket is also
called {\it special Liouville} in these coordinates.
The special Liouville coordinates were introduced
in \cite{16}, \cite{17}. The most important
case is the case of nondegenerate matrix $\eta_{ij}$.

\bt
An arbitrary Poisson bracket of form (\ref{nonl2})
is special Liouville in local coordinates
$u = (u^1,..., u^N)$ if and only if
it is Liouville with a special Liouville
function $\Phi^{ij} (u)$ and functions
$(\varphi^{\alpha})^i (u)$ of the special form
such that
\be
\eta_{ks} \Phi^{sj} (u) = {\pa F^j \over \pa u^k},
\ \ \ \ \eta_{ks} (\varphi^{\alpha})^s (u) = {\pa \psi^{\alpha}
\over \pa u^k}.
\ee
\et

In this case, for a nondegenerate matrix $(\eta_{ij})$,
we get exactly our bracket (\ref{nonloc3})
from the canonical compatible pair.

Thus our problem on compatible nonlocal Poisson
brackets of hydrodynamic type is
equivalent to to the problem of classification
of the special Liouville coordinates for nonlocal
Poisson brackets of hydrodynamic type.

\bt
An arbitrary nonlocal Poisson bracket of hydrodynamic type
of form (\ref{nonl2}) is compatible with the constant
Poisson bracket (\ref{co}) if and only if the functions
$u^i (x),$ $1 \leq i \leq N,$ and
$\eta_{ij} u^i (x) u^j (x),$ $\eta^{is} \eta_{sj} =
\delta^i_j,$ are densities of integrals in
involution with respect to the Poisson bracket (\ref{nonl2}).
\et

Note that $u^i (x),$ $1 \leq i \leq N,$ are the densities
of the annihilators of the bracket (\ref{co}),
and ${1 \over 2} \eta_{ij}
u^i (x) u^j (x)$ is the density of the momentum
of bracket (\ref{co}).

\bt
An arbitrary nonlocal Poisson bracket
of hydrodynamic type of form (\ref{nonl2})
is compatible with an arbitrary nondegenerate
local Poisson bracket of hydrodynamic type (\ref{lok})
if and only if $N$ annihilators and the momentum
of bracket (\ref{lok}) are integrals in involution with respect
to the Poisson bracket (\ref{nonl2}).
\et

\section{Integrable bi-Hamiltonian hierarchies}

Consider a pair of compatible
Hamiltonian operators of hydrodynamic type
$P^{ij}_1$ and $P^{ij}_2$, one of which,
let us assume $P^{ij}_2$, is local, and another
is an arbitrary nonlocal operator of form
(\ref{nonl2}). Apparently, one of the nonlocal
Hamiltonian operators of hydrodynamic type
always can be reduced to the canonical constant
form by series of reciprocal transformations if
it is nondegenerate (an analog of the classical
Darboux theorem in symplectic geometry), so that
the considered case is, in fact, general.
If the local Hamiltonian operator $P^{ij}_2$ is nondegenerate,
then it follows from theorem \ref{t1} that,
by local change of coordinates, the pair of
compatible Hamiltonian operators
$P^{ij}_1$ and $P^{ij}_2$ can be reduced to
the following canonical form:

\be
P^{ij}_2 [v(x)] = \eta^{ij} {d \over dx}, \label{op1}
\ee
\bea
&&
P^{ij}_1 = \left ( \eta^{is}
{\pa F^j \over \pa u^s} +
\eta^{js}
{\pa F^i \over \pa u^s} -
\eta^{is} \eta^{jk}
\sum_{\alpha = 1}^L
\varepsilon_{\alpha} {\pa \psi^{\alpha} \over \pa u^s}
{\pa \psi^{\alpha} \over \pa u^k} \right )
 {d \over dx} +  \nn\\
&&
\left ( \eta^{is} {\pa^2 F^j \over \pa u^s \pa u^k}
- \eta^{is} \eta^{jp} \sum_{\alpha = 1}^L
\varepsilon_{\alpha} {\pa^2 \psi^{\alpha} \over
\pa u^s \pa u^k} {\pa \psi^{\alpha} \over \pa u^p} \right )
\, u^k_x + \nn\\
&&
\eta^{ip} \eta^{jr}
\sum_{\alpha =1}^L
\varepsilon_{\alpha} {\pa^2 \psi^{\alpha} \over \pa u^p \pa u^k}
u^k_x \left ( {d \over dx} \right )^{-1}
{\pa^2 \psi^{\alpha} \over \pa u^r \pa u^s} u^s_x,  \label{op2}
\eea
where ($\eta^{ij}$) is an arbitrary nondegenerate
constant symmetric matrix:
$\det (\eta^{ij}) \neq 0,$ $ \eta^{ij} = {\rm const},$
$\eta^{ij} = \eta^{ji};$ $F^i (u),$  $1 \leq i \leq N,$  and
$\psi^{\alpha} (u),$ $1 \leq \alpha \leq L,$ are
smooth functions defined in a certain domain of
local coordinates such that operator (\ref{op2})
is Hamiltonian, that is, the functions
$F^i (u)$ and $\psi^{\alpha} (u)$ satisfy
the integrable equations (\ref{ass1}), (\ref{ass2})
(see theorem \ref{liu} above).

\br
{\rm It is obvious that here we can always
consider that $\eta^{ij} = \varepsilon^i \delta^{ij},$
$\varepsilon^i = 1$ for $i \leq p,$
$\varepsilon^i = - 1$ for $i > p,$
where $p$ is the positive index of inertia of the metric,
$0 \leq p \leq N$, and, in addition,
it is necessary to classify the Hamiltonian operators
(\ref{op2}) with respect to the action of
the group of motions for the corresponding
$N$-dimensional pseudo-Euclidean space ${\bf R}^N_p$,
but for our purposes it is sufficient
(and more convenient) to use the indicated above
representation for canonical compatible pair
(``conventionally canonical''
representation).}
\er

Consider the recursion operator generated
by the canonical compatible Hamiltonian operators
(\ref{op1}),
(\ref{op2}):

\bea
&&
R^i_l [v(x)] =
\left [ P_1 [v (x)] \left ( P_2 [v(x)] \right )^{-1} \right ]^i_l
=
\left (
\left ( \eta^{is}
{\pa F^j \over \pa u^s} +
\eta^{js}
{\pa F^i \over \pa u^s} - \right. \right. \nn\\
&&
\left. \eta^{is} \eta^{jk}
\sum_{\alpha = 1}^L
\varepsilon_{\alpha} {\pa \psi^{\alpha} \over \pa u^s}
{\pa \psi^{\alpha} \over \pa u^k} \right )
 {d \over dx} +
\left ( \eta^{is} {\pa^2 F^j \over \pa u^s \pa u^k}
- \eta^{is} \eta^{jp} \sum_{\alpha = 1}^L
\varepsilon_{\alpha} {\pa^2 \psi^{\alpha} \over
\pa u^s \pa u^k} {\pa \psi^{\alpha} \over \pa u^p} \right )
\, u^k_x + \nn\\
&&
\left. \eta^{ip} \eta^{jr}
\sum_{\alpha =1}^L
\varepsilon_{\alpha} {\pa^2 \psi^{\alpha} \over \pa u^p \pa u^k}
u^k_x \left ( {d \over dx} \right )^{-1}
{\pa^2 \psi^{\alpha} \over \pa u^r \pa u^s} u^s_x
\right ) \eta_{jl} \left ( {d \over dx} \right )^{-1}
\label{recur}
\eea
(what about recursion operators
generated by pairs of compatible Hamiltonian
operators, see \cite{18}--\cite{23}.

Let us apply the derived recursion operator (\ref{recur})
to the system of translations with respect to $x$,
that is, the system of hydrodynamic type
\be
u^i_t = u^i_x,  \label{transl}
\ee
which is, obviously, Hamiltonian with the
Hamiltonian operator (\ref{op1}):
\be
u^i_t =u^i_x \equiv P^{ij}_2 {\delta H \over \delta u^j (x)},
\ \ \ H = {1 \over 2} \int \eta_{jl} u^j (x) u^l (x) dx.
\ee

Any system from hierarchy
\be
u^i_{t_n} = \left ( R^n \right )^i_j u^j_x,
\ \ \ \ n \in {\bf Z}, \label{ierarkh}
\ee
is a multi-Hamiltonian integrable system.

In particular, any system of the form
\be
u^i_{t_1} = R^i_j u^j_x,
\ee
that is, the system of hydrodynamic type
\bea
&&
u^i_{t_1} =
\left (
\left ( \eta^{is}
{\pa F^j \over \pa u^s} +
\eta^{js}
{\pa F^i \over \pa u^s} - \eta^{is} \eta^{jk}
\sum_{\alpha = 1}^L
\varepsilon_{\alpha} {\pa \psi^{\alpha} \over \pa u^s}
{\pa \psi^{\alpha} \over \pa u^k} \right )
 {d \over dx} + \right. \nn\\
&&
\left ( \eta^{is} {\pa^2 F^j \over \pa u^s \pa u^k}
- \eta^{is} \eta^{jp} \sum_{\alpha = 1}^L
\varepsilon_{\alpha} {\pa^2 \psi^{\alpha} \over
\pa u^s \pa u^k} {\pa \psi^{\alpha} \over \pa u^p} \right )
\, u^k_x + \nn\\
&&
\left. \eta^{ip} \eta^{jr}
\sum_{\alpha =1}^L
\varepsilon_{\alpha} {\pa^2 \psi^{\alpha} \over \pa u^p \pa u^k}
u^k_x \left ( {d \over dx} \right )^{-1}
{\pa^2 \psi^{\alpha} \over \pa u^r \pa u^s} u^s_x
 \right ) \eta_{jl} u^l
\equiv \nn\\
&&
\left ( \eta^{is} {\pa F^j \over \pa u^s}
\eta_{jk} + {\pa F^i \over \pa u^k}
+ \eta^{is} {\pa^2 F^j \over \pa u^s \pa u^k}
\eta_{jr} u^r - \eta^{is} \sum_{\alpha = 1}^L
\varepsilon_{\alpha} {\pa \psi^{\alpha} \over \pa u^s}
{\pa \psi^{\alpha} \over \pa u^k} - \right. \nn\\
&&
\left. \eta^{is} \sum_{\alpha = 1}^L
\varepsilon_{\alpha} {\pa^2 \psi^{\alpha}
\over \pa u^s \pa u^k} \psi^{\alpha} \right )
u^k_x
\equiv
\left ( F^i (u) + \eta^{is} \eta_{jr} {\pa F^j \over \pa u^s} u^r
- \eta^{is} \sum_{\alpha = 1}^L \varepsilon_{\alpha}
{\pa \psi^{\alpha} \over
\pa u^s} \psi^{\alpha} \right )_x,  \label{canon}
\eea
where $F^i(u),$ $1 \leq i \leq N,$ and $\psi^{\alpha} (u),$
$1 \leq \alpha \leq L,$ is an arbitrary solution
of the integrable system (\ref{ass1}), (\ref{ass2}).

This system of hydrodynamic type is bi-Hamiltonian
with the pair of canonical compatible
Hamiltonian operators
(\ref{op1}),
(\ref{op2}):
\bea
&&
u^i_{t_1} =
\left (
\left ( \eta^{is}
{\pa F^j \over \pa u^s} +
\eta^{js}
{\pa F^i \over \pa u^s} -
\eta^{is} \eta^{jk}
\sum_{\alpha = 1}^L
\varepsilon_{\alpha} {\pa \psi^{\alpha} \over \pa u^s}
{\pa \psi^{\alpha} \over \pa u^k} \right )
 {d \over dx} + \right. \nn\\
&&
\left ( \eta^{is} {\pa^2 F^j \over \pa u^s \pa u^k}
- \eta^{is} \eta^{jp} \sum_{\alpha = 1}^L
\varepsilon_{\alpha} {\pa^2 \psi^{\alpha} \over
\pa u^s \pa u^k} {\pa \psi^{\alpha} \over \pa u^p} \right )
\, u^k_x + \nn\\
&&
\left. \eta^{ip} \eta^{jr}
\sum_{\alpha =1}^L
\varepsilon_{\alpha} {\pa^2 \psi^{\alpha} \over \pa u^p \pa u^k}
u^k_x \left ( {d \over dx} \right )^{-1}
{\pa^2 \psi^{\alpha} \over \pa u^r \pa u^s} u^s_x
\right )
{\delta H_1 \over \delta u^j (x)}, \nn\\
&&
H_1 = {1 \over 2} \int \eta_{jl} u^j (x) u^l (x) dx, \label{eq1}
\eea
\be
u^i_{t_1} = \eta^{ij} {d \over dx}
{\delta H_2 \over \delta u^j (x)}, \ \ \
H_2 =  \int \left (\eta_{jk} F^k (u (x)) u^j (x)
- {1 \over 2} \sum_{\alpha = 1}^L
\varepsilon_{\alpha} (\psi^{\alpha} (u))^2 \right )
dx. \label{eq2}
\ee

\medskip

\begin{flushleft}
Centre for Nonlinear Studies,\\
L.D.Landau Institute for Theoretical Physics, \\
Russian Academy of Sciences\\
e-mail: mokhov@mi.ras.ru; mokhov@landau.ac.ru\\
\end{flushleft}

\end{document}